\documentclass[twoside]{article} 

\usepackage{jmlr2e} 
\usepackage{graphicx} 
\usepackage{subfigure}

\usepackage{natbib}

\usepackage{algorithm}
\usepackage{algorithmic}

\usepackage{hyperref}

\usepackage{amsmath}
\usepackage{pgfplots}
\usepackage{color}
\usepackage{needspace}

\usepackage{flushend}
\usepackage{multirow}
\usepackage{rotating}
\usepackage{appendix}
\usepackage{tikz,pgfplots}

\usepackage{standalone} 

\newif\ifcompiletikz
\compiletikztrue

\newtheorem{assumption}{Assumption}

\input{mysymbol.sty}

\def\whp{\text{w.h.p}}
\newcommand{\QED}{\hfill\ensuremath{\blacksquare}}

\begin{document}




\ShortHeadings{Large Scale ERM via Truncated Adaptive Newton Method}{Eisen, Mokhtari, and Ribeiro}
\firstpageno{1}



\title{Large Scale Empirical Risk Minimization via Truncated Adaptive Newton Method}

\author{\name Mark Eisen \email maeisen@seas.upenn.edu \\
         \name Aryan Mokhtari  \email aryanm@seas.upenn.edu \\
       \name Alejandro Ribeiro \email aribeiro@seas.upenn.edu  \\
       \addr Department of Electrical and Systems Engineering\\
       University of Pennsylvania\\
       Philadelphia, PA 19104, USA}

\editor{}

\maketitle


\begin{abstract}
We consider large scale empirical risk minimization (ERM) problems, where both the problem dimension and variable size is large. In these cases, most second order methods are infeasible due to the high cost in both computing the Hessian over all samples and computing its inverse in high dimensions. In this paper, we propose a novel adaptive sample size second-order method, which reduces the cost of computing the Hessian by solving a sequence of ERM problems corresponding to a subset of samples and lowers the cost of computing the Hessian inverse using a truncated eigenvalue decomposition. We show that while we geometrically increase the size of the training set at each stage, a single iteration of the truncated Newton method is sufficient to solve the new ERM within its statistical accuracy. Moreover, for a large number of samples we are allowed to double the size of the training set at each stage, and the proposed method subsequently reaches the statistical accuracy of the full training set approximately after two effective passes.
 In addition to this theoretical result, we show empirically on a number of well known data sets that the proposed truncated adaptive sample size algorithm outperforms stochastic alternatives for solving ERM problems.

\end{abstract}


\section{Introduction}

The recent advances in large scale machine learning focus largely on solving of the expected risk minimization problem, in which a set of model parameters of dimension $p$ are found that minimize an expected loss function. More typically, the expected loss function is taken with respect to an unknown probability distribution. Therefore, the expected risk is estimated with a statistical average over $N$ samples, where $N$ is very large. The minimization over the statistical average is called the empirical risk minimization (ERM) problem. Indeed the computational complexity of solving ERM problems depends on the size of the dataset $N$ if we use deterministic methods that operate on the full dataset at each iteration. Stochastic optimization has long been used to reduce this cost. Among the stochastic optimization methods are stochastic gradient descent algorithm\citep{robbins1951stochastic,bottou2010large}, Nesterov-based methods \citep{nesterov2007gradient, beck2009fast}, variance reduction methods \citep{johnson2013accelerating,nguyen2017sarah}, stochastic average gradient algorithms \citep{roux2012stochastic,defazio2014saga}, stochastic majorization-minimization algorithms \citep{defazio2014finito,mairal2013stochastic}, hybrid methods \citep{konecny2013semi}, and dual coordinate methods \citep{ shalev2013stochastic, shalev2016accelerated}. Although these stochastic first-order methods succeed in reducing the computational complexity of deterministic methods, they suffer from slow convergence in ill conditioned problems. This drawback inspires the development of stochastic second order methods, which improve upon the performance of first-order methods by using a curvature correction. These include subsampled Hessian \citep{erdogdu2015convergence,pilanci2015newton,roosta2016sub,roosta2016sub2}, incremental Hessian updates \citep{gurbuzbalaban2015globally},  stochastic dual Newton ascent \citep{DBLP:conf/icml/QuRTF16}, and stochastic quasi-Newton methods \citep{schraudolph2007stochastic,mokhtari2014res,JMLR:v16:mokhtari15a, lucchi2015,moritz2015linearly, byrd2016stochastic, mokhtari2017iqn}. However, many of these cannot improve upon the asymptotic convergence rate of linearly convergent first-order methods. Those that do improve rates do so at great computational cost due to the computation of the full gradient and the Hessian inverse, both computationally prohibitive when $N$ and $p$ are large, as is often true in real datasets.

Each of these methods solve the ERM problem using the full sample set. However, the samples in ERM problems are drawn from a common distribution, so a smaller ERM problem using a smaller subset of samples should have a solution close to that of the full problem. Solving a sequence of smaller ERM problems using a subset of samples thus may reduce computational complexity. The work in \citep{mokhtari2016adaptive}, for instance, reduces complexity of Newton's method by adaptively increasing sample size, but remains impractical for high dimensional problems due to the inverse computation. However, a key insight in ERM problems is that only the empirical average is optimized, and not the true expected loss function. Indeed this error between the optimal empirical loss and expected loss, called the \emph{statistical accuracy}, can be made arbitrarily small by increasing the number of samples $N$. More importantly, however, is that the unavoidable presence of such an error in fact grants us latitude in the accuracy of our optimization methods. Given that we can only optimize the expected loss function up to a certain accuracy, we can employ further approximation techniques with negligible loss, so long as the error induced by these approximations is small relative to that of the stochastic approximation. It is therefore possible to significantly reduce the computational cost of our optimization methods without negatively impacting the overall accuracy, thus making second order methods a viable option.

%
 
%

In this paper, we propose a novel adaptive sample size second-order method, which reduces the cost of computing the Hessian by solving a sequence of ERM problems corresponding to a subset of samples and lowers the cost of computing the inverse of the Hessian by using a truncated eigenvalue decomposition. In the presented scheme, we increase the size of the training set geometrically, while solving subproblems up to their statical accuracy. We show that we can increase the sample size in such a manner that a single iteration of a truncated Newton method is sufficient to solve the increased sample size ERM problem to within its statistical accuracy. This is achieved by using the quadratic convergence region of the new ERM problem. While the proposed method does not converge at a purely quadratic rate due to truncation of the Newton step, the additional linear term incurred by the approximation can be made negligible with respect to the statistical accuracy of the data set. The resulting $k$-Truncated Adaptive Newton ($k$-TAN) method uses a rank-$k$ approximation of the Hessian to solve high dimensional ERM problems with large sample sizes up to the statistical accuracy at a significantly lower overall cost than existing optimization techniques. Specifically, we demonstrate that, in many cases, we can reach the statistical accuracy of the problem with a total computational cost of $\ccalO((2N + (\log_2N) \log k) p^2 )$ using rank $k$ approximations of the Hessian inverse.

\section{Problem Formulation}

We consider in this paper the empirical risk minimization (ERM) problem for a convex function $f(\bbx, z)$, where $z$ is a realization of a random variable $Z$. More specifically, we seek the optimal variable $\bbx \in \reals^p$ that minimizes the expected loss $L(\bbx) := \bbE_Z [ f(\bbx,z)]$. Define $\bbx^*$ as the variable that minimizes the expected loss, i.e.
\begin{equation} \label{eq_orig_problem}
\bbx^* := \argmin_{\bbx} \bbE_Z [ f(\bbx,z)].
\end{equation}
In general, the problem in \eqref{eq_orig_problem} cannot be solved without knowing the distribution of $Z$. As an alternative, we traditionally consider the case that we have access to $N$ samples of $Z$, labelled $z_1, z_2, \hdots, z_n$. Define then the functions $f_i(\bbx) = f(\bbx, z_i)$ for $i = 1, 2, \hdots, N$ and an associated empirical risk function $L_n := (1/n) \sum_{i=1}^n f_i(\bbx)$ as the statistical mean over the first $n \leq N$ samples. We say that function $L_n(\bbx)$ approximates the original expected loss $L(\bbx)$ with statistical accuracy $V_n$ if the difference between the empirical risk function $L_n(\bbx)$ and the expected loss $L(\bbx)$ is upper bounded by $V_n$ for all $\bbx$ with high probability (w.h.p.). The statistical accuracy $V_n$ is typically bounded by $V_n = \ccalO(1/\sqrt{n})$ \citep{vapnik1998statistical} or the stronger $V_n = \ccalO(1/n)$ for a set of common problems \citep{bartlett2006convexity,bousquet2008tradeoffs}. 

Observe that the sampled loss function $L_n$ is of an order $V_n$ difference from the true loss function $L$ and, consequently, any additional change of the same order has negligible effect. It is therefore common to regularize non-strongly convex loss functions $L_n$ by a term of order $V_n$. We then seek the minimum argument of the regularized risk function $R_n$,
\begin{equation} \label{eq_reg_problem}
\bbx_n^* := \argmin_{\bbx} R_n(\bbx) := \argmin_{\bbx} \frac{1}{n} \sum_{i=1}^n f_i(\bbx) + \frac{cV_n}{2}\|\bbx\|^2 ,
\end{equation}
where $c$ is a scalar constant. The solution $\bbx_n^*$ minimizes the regularized risk function using the first $n$ samples, which is of order $V_n$ from the expected loss function $L$. It follows then that by setting $n=N$ we find a solution $\bbx_N^*$ in \eqref{eq_reg_problem} that solves the original problem in \eqref{eq_orig_problem} up to the statistical accuracy of using all $N$ samples.

The problem in \eqref{eq_reg_problem} is strongly convex can be solved using any descent method. In particular, Newton's method uses a curvature-corrected gradient to iteratively update a variable $\bbx$, and is known to converge to the optimal argument $\bbx_n^*$ at a very fast quadratic rate. To implement Newton's method, it is necessary to compute the gradient $\nabla R_n(\bbx)$ and Hessian $\nabla^2 R_n(\bbx)$ as
\begin{align}\label{eq_grad_hessian}
\nabla R_n(\bbx) =  \frac{1}{n} \sum_{i=1}^n \nabla f_i(\bbx) + cV_n\bbx, \qquad  \nabla^2 R_n(\bbx)= \frac{1}{n} \sum_{i=1}^n \nabla^2 f_i(\bbx) + cV_n \bbI.
\end{align}
The variable $\bbx$ is updated in Newton's method as
\begin{align}\label{eq_newton_update}
\bbx^+ = \bbx -   \nabla^{-2} R_n(\bbx) \nabla R_n(\bbx).
\end{align}
Solving \eqref{eq_orig_problem} to the full statistical accuracy $V_N$ (i.e. solving \eqref{eq_reg_problem} for $n=N$) using Newton's method would then require the computation of individual gradients and Hessians for $N$ functions $f_i$ for computational cost of $\ccalO(N p^2)$ at each iteration. Furthermore, the computation of the Hessian inverse in \eqref{eq_newton_update} requires a cost of $\ccalO(p^3)$, bringing at total of $\ccalO(Np^2 + p^3)$ for an iteration of Newton's method using the whole dataset. For large $n$ and $p$, this may become computationally infeasible. In this paper we show how this complexity can be reduced by gradually increasing the sample size $n$ and approximating the inverse of the respective Hessian $\nabla^2 R_n(\bbx)$.

\section{$k$-Truncated Adaptive Newton ($k$-TAN) Method}\label{sec_ada_newton}

We propose the $k$-Truncated Adaptive Newton ($k$-TAN) as a low cost alternative to solving \eqref{eq_orig_problem} to its statistical accuracy. In the $k$-TAN method, at each iteration we start from a point $\bbx_m$  within the statistical accuracy of $R_m$, i.e. $R_m (\bbx_m)-R_m (\bbx_m^*) \leq V_m$. We geometrically increase the sample size to $n = \alpha m$, where $\alpha>1$, and compute $\bbx_n$ using an approximated Newton method on the increased sample size risk function $R_n$. More specifically, we update a decision variable $\bbx_m$ associated with $R_m$ to a new decision variable $\bbx_n$ associated with $R_n$ with the Newton-type update
\begin{align} \label{eq_update}
\bbx_n = \bbx_m - \hbH_{n,k}^{-1} \nabla R_n(\bbx_m),
\end{align}
where $\hbH_{n,k}$ is a matrix approximating the Hessian $\nabla^2 R_n(\bbx_m)$ and parametrized by $k \in \{1,2,\hdots,p\}$. In particular we are interested in an approximation $\hbH_{n,k}$ whose inverse $\hbH^{-1}_{n,k}$ can be computed with complexity less than $\ccalO(p^3)$.  To define such a matrix, consider $\mu_1 \geq \mu_2 \geq \hdots \geq \mu_p$ to be the eigenvalues of the Hessian of empirical risk $\nabla^2 L_n (\bbx_m)$, with associated eigenvectors $\bbv_1, \bbv_2, \hdots, \bbv_p$. We perform an eigenvalue decomposition of $\nabla^2  L_n(\bbx_m) = \bbU \bbSigma \bbU^T$, where $\bbU := [\bbv_1, \hdots, \bbv_p]\in \reals^{p \times p}$ and $\bbSigma := \diag(\mu_1, \hdots, \mu_p)\in \reals^{p \times p}$. We can then define the truncated eigenvalue decomposition with rank $k$ as $\hat{\nabla}^2 L_n (\bbx_m) := \bbU_k \bbSigma_k \bbU_k^T$, where $\bbU_k := [\bbv_1, \hdots, \bbv_k] \in \reals^{p \times k}$ and $\bbSigma_k := \diag(\mu_1, \hdots, \mu_k) \in \reals^{k \times k}$. The full approximated Hessian $\hbH_{n,k}$ is subsequently defined as the rank $k$ approximation of $\nabla^2  L_n(\bbx_m)$ regularized by $c V_n \bbI$, i.e.
\begin{align}
\hbH_{n,k} &:= \bbU_k \bbSigma_k \bbU_k^T + cV_n \bbI \label{eq_a_hessian}.
\end{align}
The inverse of the approximated Hessian $\hbH_{n,k}$ can then be computed directly using $\bbU_k$ and $\bbSigma_k$ as
\begin{align}
\hbH^{-1}_{n,k} &:= \bbU_k [(\bbSigma_k+cV_n \bbI)^{-1} - (cV_n\bbI)^{-1}]\bbU_k^T + (cV_n)^{-1} \bbI. \label{eq_a_hessian_inverse}
\end{align}

Observe that setting $k=p$, i.e., full Hessian inverse, recovers the AdaNewton method in \citep{mokhtari2016adaptive}. To understand how we may determine $k$, consider that the full Hessian computed in \eqref{eq_grad_hessian} is $\nabla^2  L_n(\bbx_m)$ regularized by $ c V_n \bbI$. Therefore, the eigenvalues of $\nabla^2  L_n(\bbx_m)$ less than $c V_n$ are made negligible by the regularization, and can be left out of the approximation. We thus select the $k$ largest eigenvalues of the Hessian which are larger than $\rho c V_n$ for some truncation parameter $0 < \rho < 1$.

To analyze the computational complexity of \eqref{eq_a_hessian_inverse}, observe that the inverse computation in \eqref{eq_a_hessian_inverse} requires only the inversion of diagonal matrices, and thus the primary cost in computing the $k$ largest eigenvalues $\bbSigma_k$ and associated eigenvectors $\bbU_k$. Indeed, the truncated eigenvalue decomposition $\{ \bbU_k, \bbSigma_k\}$ can in general be computed with at most complexity $\ccalO(k p^2)$, with recent randomized algorithms even finding $\{ \bbU_k, \bbSigma_k\}$ with complexity $\ccalO(p^2 \log k)$ \citep{halko2011finding}. This results in a total cost of, at worst, $\ccalO((\log k+n)p^2)$ to perform the update in \eqref{eq_update}, thus removing a $\ccalO(p^3)$ cost.

In this paper we aim to show that while we geometrically increase the size of the training set, a single iteration of the truncated Newton method in \eqref{eq_update} is sufficient to solve the new risk function within its statistical accuracy. To state this result we first assume the following assumptions hold. 

\begin{assumption}\label{ass_convexity}
The loss functions $f(\bbx,\bbz)$ are convex with respect to $\bbx$ for all values of $\bbz$. Moreover, their gradients $\nabla f(\bbx,\bbz)$ are Lipschitz continuous with constant $M$.
\vspace{-3mm}
\end{assumption}

\begin{assumption}\label{ass_self_concor}
The loss functions $f(\bbx,\bbz)$ are self-concordant with respect to $\bbx$ for all $\bbz$. 
\vspace{-3mm}
\end{assumption}

\begin{assumption}\label{ass_grad_cond}
The difference between the gradients of the empirical loss $L_n$ and the statistical average loss $L$ is bounded by $V_n^{1/2}$ for all $\bbx$ with high probability,
\begin{align}\label{eqn_loss_minus_erm_2}
   \sup_{\bbx}\|\nabla L(\bbx) - \nabla L_{n}(\bbx) \|  \leq V_n^{1/2},  \qquad\whp.
\end{align}
\vspace{-5mm}
\end{assumption}

Based on Assumption \ref{ass_convexity}, we obtain that the regularized empirical risk gradients $\nabla R_n$ are Lipschitz continuous with constant $M+cV_n$. 
%
 Assumption \ref{ass_self_concor} states the loss functions are additionally self concordant which is a customary assumption in the analysis of second-order methods. It also follows that the functions $R_n$ are therefore self concordant. Assumption \ref{ass_grad_cond} bounds the difference between gradients of the expected loss and the empirical loss with $n$ samples by $V_n^{1/2}$. This is a reasonable bound for the convergence of gradients to their statistical averages using the law of large numbers.



 We are interested in establishing the result that, as we increase $n$ at each step, the $k$-TAN method stays in the quadratic region of the the associated risk function. More explicitly, we wish to show the sample size can be increased from $m$ to $n = \alpha m$ such that $\bbx_m$ is in the quadratic region of $R_n$. Moreover, if $\bbx_m$ is indeed in the quadratic region of $R_n$, then we demonstrate that a single step of $k$-TAN as in \eqref{eq_update} produces a point $\bbx_n$ that is within the statistical accuracy $V_n$ of the risk $R_n$. 

%
\begin{theorem}\label{theorem_main_result}
Consider the $k$-TAN method defined in \eqref{eq_update}-\eqref{eq_a_hessian_inverse} and suppose that the constant $k$ for low rank factorization is defined as $k = \min \{k\ |\ \mu_{k+1} \leq \rho cV_n\}$ where  $\rho$ is a free parameter chosen from the interval $(0,1]$. Further consider the variable $\bbx_m$ as a $V_m$-optimal solution of the risk $R_{m}$, i.e., a solution such that $R_{m}(\bbx_m)- R_{m}(\bbx_m^*) \leq V_m$. Let $n=\alpha m > m$ and suppose Assumptions \eqref{ass_convexity}-\eqref{ass_grad_cond} hold. If the sample size $n$ is chosen such that the following conditions
\begin{equation}\label{cond_1}
   \left(\frac{2(M+cV_m)V_{m}}{cV_n}\right)^{1/2}\!
                     + \frac{2(n-m)}{nc^{1/2}}     
                     + \frac{\left((2+\sqrt{2})c^{1/2} +c\|\bbx^*\|\right)(V_m-V_n)}{(cV_n)^{1/2}}\leq \frac{1}{4}
\end{equation}
\begin{equation}\label{cond_2}
\frac{16}{(3 - \rho)^4}\left[ 36 K^2(1+\rho)^2 V_m^2 + 30 K^{3/2} \rho(1+\rho)V_m^{3/2} + 6K\rho^2 V_m\right]
\leq V_n
\end{equation}
are satisfied, where $K = 3 + (2+c \| \bbx^*\|^2/2)(1 - 1/\alpha)$, then the variable $\bbx_n$ computed from \eqref{eq_update} has the suboptimality of $V_n$ with high probability, i.e., 
\begin{equation}\label{imp_result}
 R_{n}(\bbx_n)- R_{n}(\bbx_n^*) \leq V_n, \qquad \whp.
\end{equation}
 \end{theorem}

%

The result in Theorem \ref{theorem_main_result} establishes the required conditions to guarantee that the iterates $\bbx_n$ always stay within the statistical accuracy of the risk $R_n$. The expression in \eqref{cond_1} provides a condition on growth rate $\alpha = n/m$ to ensures that iterate $\bbx_m$, which is a $V_m$-suboptimal solution for $R_m$, is within the quadratic convergence neighborhood of Newton's method for $R_n$. The second condition in \eqref{cond_2} ensures that a single iteration of $k$-TAN is sufficient for the updated variable $\bbx_n$ to be within the statistical accuracy of $R_n$. Note that the first term in the left hand side of \eqref{cond_2} is quadratic with respect to $V_m$ and comes from the quadratic convergence of Newton's method, while the second and third terms of respective orders $V_m^{3/2}$ and $V_m$ are the outcome of Hessian approximation. Indeed, these terms depend on $\rho$, which is the upper bound on ratio of the discarded eigenvalues $\mu_{k+1},\hdots,\mu_p$ to the regularization $cV_n$. The truncation must be enough such that $\rho$ is sufficiently small to make \eqref{cond_2} hold. It is worth mentioning, as a sanity check, if we set $\rho=0$ then we will keep all the eigenvalues and recover the update of Newton's method which makes the non-quadratic terms in \eqref{cond_2} zero.

The conditions in Theorem \ref{theorem_main_result} are cumbersome but can be simplified if we focus on large $m$ and assume that the inequality $V_m\leq \alpha V_n$ holds for $n=\alpha m$. Then, \eqref{cond_1} and \eqref{cond_2} can be simplified to 
\begin{equation}\label{eq_cond_simple}
\left( \frac{2\alpha M}{c} \right)^{\!1/2}\!\! +\frac{2(\alpha-1)}{\alpha c^{1/2}}\leq \frac{1}{4},\quad  \text{and} \quad  \frac{96 [3 + (2+c \| \bbx^*\|^2/2)(1 - 1/\alpha)]\rho^2}{(3-\rho)^2}\leq \frac{1}{\alpha},
\end{equation}
respectively. 
Observe that first condition is dependent of $\alpha$ and the second condition depends on $\alpha$ and $\rho$. Thus, a pair $(\alpha,\rho)$ must be chosen that satisfies \eqref{eq_cond_simple} for the result in Theorem \ref{theorem_main_result} to hold. We point out one such pair as the parameters $\alpha=2$, $c>16(2\sqrt{M} + 1)^2$ and $\rho = 9/(21\sqrt{c \|\bbx^*\|^2+16}+3)$. Consequently, when $m$ is large we may double the sample size with each update in until $n=N$, after which we will have obtained a point $\bbx_N$ such that $R_N(\bbx_N) -R_N(\bbx_N^*) \leq V_N$. After $\log_{2} N$ iterations (roughly $2N$ samples processed), we solve the full risk function $R_N$ to within the statistical accuracy $V_N$. At each iteration, the truncated inverse step requires cost $\ccalO(p^2\log k)$. Computing Hessians over $2N$ samples requires cos $\ccalO(2Np^2)$, resulting in a total complexity of $\ccalO(p^2(2N + \log_2 N \log k)$.

In practice, these may be chosen in a backtracking manner, in which the iterate $\bbx_m$ is updated using an estimate $(\alpha,\rho)$ pair. If the resulting iterate $\bbx_n$ is not in statistical accuracy $R_N(\bbx_N) -R_N(\bbx_N^*) \leq V_N$, the increase factor $\alpha$ is decreased by factor $\beta$ and $\rho$ is decreased by factor $\delta$. Since $R_N(\bbx_N^*)$ is not known in practice, the suboptimality can be upper bounded using strong convexity as $R_N(\bbx_N) -R_N(\bbx_N^*)\leq \| \nabla R_n(\bbx_n) \|^2/ (2 cV_n)$.

The resulting method is presented in Algorithm \ref{alg:AdaNewton}. After preliminaries and initializations in Steps 1-4, the backtracking loop starts in Step 6 with the sample size increase by rate $\alpha$. After computing the gradient and Hessian in Step 7 and 8, the low rank decomposition with rate $k =\min \{k|\mu_{k+1} \leq \rho cV_n\}$ in Step 9. The $k$-TAN update is then performed with \eqref{eq_a_hessian_inverse}-\eqref{eq_update} in Steps 9 and 10. The factors $\alpha$ and $\rho$ are then decreased using the backtracking parameters and the statistical accuracy condition is checked. We stress that, while $\nabla R_n(\bbx_n)$ must be computed to check the exit condition in Step 13, the gradient for these samples must be computed in any case in the following iteration, so no additional computation is added by this step.

%
{\linespread{1.3}
\begin{algorithm}[t] \begin{algorithmic}[1]
\STATE \textbf{Parameters:} Sample size increase constants $\alpha_0>1$, $\rho_0 < 1$  and $0<\beta,\delta<1$
\STATE \textbf{Input:} Initial sample size $n=m_0$ and 
                       argument $\bbx_{n} = \bbx_{m_0}$ with 
                       $\| \nabla R_{n}(\bbx_n)\| < (\sqrt{2 c}) V_n$
\WHILE [main loop]{$n\leq N$} 
   \STATE Update argument and index:\ $\bbx_m=\bbx_n$ and $m=n$. 
          Reset factor $\alpha=\alpha_0$, $\rho = \rho_0$ .     
   \REPEAT  [sample size backtracking loop] 
      \STATE Increase sample size:\ $n=\min\{\alpha m, N\}$. 
      \STATE Compute gradient [cf. \eqref{eq_grad_hessian}]: \ 
             ${
                 \nabla R_n (\bbx_m) = (1/n)\sum_{ki=1}^{n}\nabla f(\bbx_m,z_i) + cV_n \bbx_m
             }$
      \STATE Compute Hessian: \
             ${
                 \nabla^2 L_n = (1/n)\sum_{i=1}^{n} \nabla^2 f (\bbx_m,z_i)
             }$
       \STATE Find low rank decomposition \citep{halko2011finding}:
                    ${
                 \ \hat{\nabla}^2 L_n = \bbU_k \bbSigma_k \bbU_k^T \ \text{ for } k = \min \{k|\mu_{k+1} \leq \rho cV_n\}}$
        \STATE Compute $  \hbH^{-1}_{n,k}$ [cf. \eqref{eq_a_hessian_inverse}]: 
             ${
                 \hbH^{-1}_{n,k} = \bbU_k [(\bbSigma_k+cV_n \bbI)^{-1} - (cV_n)^{-1}\bbI]\bbU_k^T + (cV_n)^{-1} \bbI
             }$
      \STATE Newton Update [cf. \eqref{eq_update}]: \
             ${
                 \bbx_n = \bbx_m - \hbH_{n,k}^{-1} \nabla R_n (\bbx_m)
             }$
      \STATE Backtrack sample size increase $\alpha=\beta\alpha$, truncation factor $\rho = \delta \rho$.       
   \UNTIL {$\| \nabla R_{n}(\bbx_n)\| < (\sqrt{2 c}) V_n$} 
\ENDWHILE
\end{algorithmic}
\caption{{$k$-TAN}}\label{alg:AdaNewton} \end{algorithm}}

\section{Convergence Analysis}

We study the convergence properties of the $k$-TAN method and in particular prove the result in Theorem 1. Namely, we show that the update described in \eqref{eq_update} produces a variable $\bbx_n$ that is within the statistical accuracy $V_n$. This, in turn, implies all future updates will be within the statistical accuracy of their respective increased sample functions $R_n$ until the full set of $N$ samples is used, at which point we will have reached a point $\bbx_N$ that is within the statistical accuracy $V_N$ of problem \eqref{eq_orig_problem}.

\subsection{Preliminaries}

Before proceeding with the analysis of $k$-TAN, we first present two propositions that relate current iterate $\bbx_m$ to the suboptimality and quadratic convergence region to the increased sample size risk $R_n$. Define $S_n(\bbx)$ to be the $n$-suboptimality of point $\bbx$ with respect to $R_n$, i.e.
\begin{align}
S_n(\bbx) := R_n(\bbx) - R_n(\bbx_n^*),
\end{align}
where $\bbx_n^*$ is the point that minimizes $R_n$. We establish in the following proposition a bound on the $n$-suboptimality of $\bbx_m$ from the difference in sample sizes $m$ and $n$ and their associated statistical accuracies. The proof can be found in \citep[Proposition 3]{mokhtari2016adaptive}.

%
\begin{proposition}\label{prop_sub_m}
Consider a point $\bbx_m$ that minimizes the risk function $R_m$ to within its statistical accuracy $V_m$, i.e. $S_m(\bbx_m) \leq V_m$. If the sample size is increased from $m$ to $n=\alpha m$ and $V_m = \alpha V_n$, the sub-optimality $S_n(\bbx_m) = R_n(\bbx_m) - R_n(\bbx_n^*)$ is upper bounded as
%
\begin{align} \label{eq_prop_sub_m_2}
S_n(\bbx_m) \leq K V_m := \left[ 3 + \left(2+\frac{c}{2} \| \bbx^*\|^2\right)\left(1 - \frac{1}{\alpha} \right) \right] V_m  \quad \text{w.h.p.}
\end{align}
%
%
\end{proposition}
Proposition \ref{prop_sub_m} demonstrates a bound on the $n$-suboptimality $S_n(\bbx_m)$ of a point $\bbx_m$ whose $m$-suboptimality is within statistical accuracy $V_m$. It is also necessary to establish conditions on increase rate $\alpha$ such the $\bbx_m$ is also in the quadratic convergence region of $R_n$. Traditional analysis of Newton's method characterizes quadratic convergence in terms of the Newton decrement $\lambda_n (\bbx) := \|\nabla^2 R_n(\bbx)^{-1/2}  \nabla R_n(\bbx) \|$. The iterate $\bbx$ is said to be in the quadratic convergence region of $R_n$ when $\lambda_n(\bbx) < 1/4$---see \citep[Chapter 9.6.4]{boyd04}. The conditions for current iterate $\bbx_m$ to be within this region are presented in the following proposition. The proof can be found in \citep[Proposition 5]{mokhtari2016adaptive}.

%
\begin{proposition}\label{main_prop}
Define $\bbx_{m}$ as an $V_m$ optimal solution of the risk $R_{m}$, i.e., ${R_{m}(\bbx_m)-R_{m}(\bbx_m^*) \leq V_m}$. In addition, define $\lambda_n (\bbx):=\left(\nabla R_{n}(\bbx) ^T \nabla^2 R_{n}(\bbx) ^{-1} \nabla R_{n}(\bbx) \right)^{1/2}$ as the Newton decrement of variable $\bbx$ associated with the risk $R_{n}$. If Assumption \ref{ass_convexity}-\ref{ass_grad_cond} hold, then Newton's method at point $\bbx_m$ is in the quadratic convergence phase for the objective function $R_{n}$, i.e., $\lambda_n(\bbw_m)<1/4$, if we have
\begin{equation}\label{prop_result}
   \left[\frac{2(M+cV_m)V_{m}}{cV_n}\right]^{1/2}\!
                           +\frac{{2(n-m)}  }{nc^{1/2}} + \frac{(\sqrt{2c}+2\sqrt{c}+c\|\bbw^*\|)(V_m-V_n)}{(cV_n)^{1/2}}\leq \frac{1}{4}
\quad \text{w.h.p.}
\end{equation}
\end{proposition}

\subsection{Analysis of $k$-TAN}

To analyze the $k$-TAN method, it is necessary to study the error incurred from approximating the Hessian inverse in \eqref{eq_a_hessian} with rank $k$. Because we are only interested in solving each risk function $R_n$ to within its statistical accuracy $V_n$, however, some approximation error can be afforded. In the following Lemma, we characterize the error between an approximate and exact Newton steps using the chosen rank $k$ of the approximation and the associated eigenvalues of the Hessian.
%
%
\begin{lemma}\label{lemma_sub_bound}
Consider the $k$-TAN update in \eqref{eq_update}-\eqref{eq_a_hessian_inverse} for some $k = \{0,1,\hdots,p\}$. Define $\epsilon_n := \mu_{k+1}/(cV_n)$. The norm of difference in the $k$-TAN step $\hbH_{n,k}^{-1}\nabla R_n(\bbx_m)$ and the exact Newton step $\bbH_n^{-1}\nabla R_n(\bbx_m)$---where $\bbH_n = \hbH_{n,p} =  \nabla^2 R_{n}(\bbx_m)$---can  be upper bounded as
\begin{align}\label{eq_lemma_sub_bound}
\| \hbH_{n,k}^{-1}\nabla R_n(\bbx_m) - \bbH_n^{-1}\nabla R_n(\bbx_m) \| \leq \epsilon_n \| \bbH_n^{-1}\nabla R_n(\bbx_m)\|.
\end{align}
\end{lemma}

The result in Lemma \ref{lemma_sub_bound} gives us an upper bound on the error incurred in single iteration of a rank $k$ approximation of the Newton step versus an exact Newton step. To make $\epsilon_n$ small, a sufficiently large $k$ must be chosen such that $\mu_{k+1}$ is in the order of $V_n$. The size of $k$ will therefore depend on the distribution of the eigenvalues of particular empirical risk function. However, in practical datasets of high dimension, it is often the case that most eigenvalues of the Hessian will be close to 0, in which case $k$ can be made very small. This trend is supported by our numerical experiments on real world data sets in Section \ref{sec_numerical_results} and the Appendix of this paper.

  With the results of Proposition \ref{prop_sub_m} and Lemma \ref{lemma_sub_bound} in mind, we can characterize the $n$-suboptimality of the updated variable $\bbx_n$ from \eqref{eq_update}. This is stated formally in the following Lemma.

\begin{lemma}\label{lemma_sub_n}
 Consider the $k$-TAN update in \eqref{eq_update}-\eqref{eq_a_hessian_inverse}. If $\bbx_m$ is in the quadratic neighborhood of $\bbR_n$, i.e. $\lambda_n (\bbx_m) < 1/4$, then the $n$-suboptimality $S_n(\bbx_n) = R_n(\bbx_n) - R(\bbx_n^*)$ can be upper bounded by
\begin{align} \label{eq_lemma_sub_n}
S_n(\bbx_n) \leq \frac{16}{(3 - \epsilon_n)^4}\left[ 36(1+\epsilon_n)^2 S_n(\bbx_m)^2 + 30 \epsilon_n(1+\epsilon_n)S_n(\bbx_m)^{3/2} + 6\epsilon_n^2 S_n(\bbx_m) \right].
\end{align}
\end{lemma}

With Lemma \ref{lemma_sub_n} we establish a bound on $n$-suboptimality of the $\bbx_n$ obtained from the $k$-TAN update in \eqref{eq_update}. Observe that this bounded by terms proportional to the $n$-suboptimality of the previous point, $S_n(\bbx_m)$. We can then establish that $S_n(\bbx_n)$ is indeed upper bounded by the statistical accuracy $V_n$ by combing the results in \eqref{eq_prop_sub_m_2} and \eqref{eq_lemma_sub_n} to obtain Theorem \ref{theorem_main_result}. The proofs of Theorem \ref{theorem_main_result} and all supporting Lemmata are provided in the Appendix.

\section{Experiments}\label{sec_numerical_results}

We compare the performance of the $k$-TAN method to existing optimization methods on  large scale machine learning problems of practical interest. In particular, we consider a regularized logistic loss function, with regularization parameters $V_n = 1/n$ and $c=1$. The $k$-TAN method is compared against the second order method AdaNewton \citep{mokhtari2016adaptive} and two first order methods---SGD and SAGA \citep{defazio2014saga}. Here, we study the performance of these methods on the logistic regression problem for two datasets. First, the GISETTE handwritten digit classification from the NIPS 2003 feature selection challenge and, second, the well-known RCV1 dataset for classifying news stories from the Reuters database. We perform additional experiments on the ORANGE dataset for KDD Cup 2009 and the BIO dataset for KDD Cup 2004.

 \begin{figure}[t]
	\centering
          \includegraphics[width=0.45\textwidth]{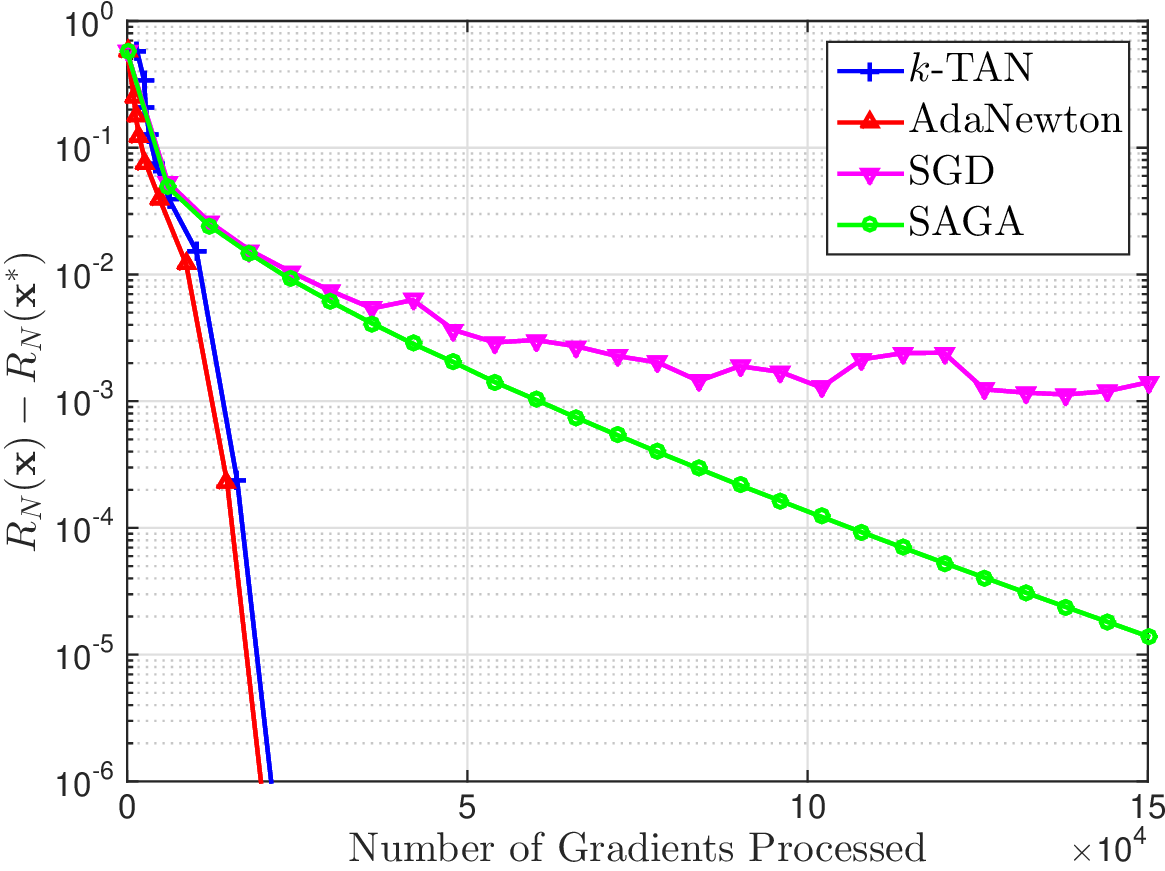}\qquad 
            \includegraphics[width=0.45\textwidth]{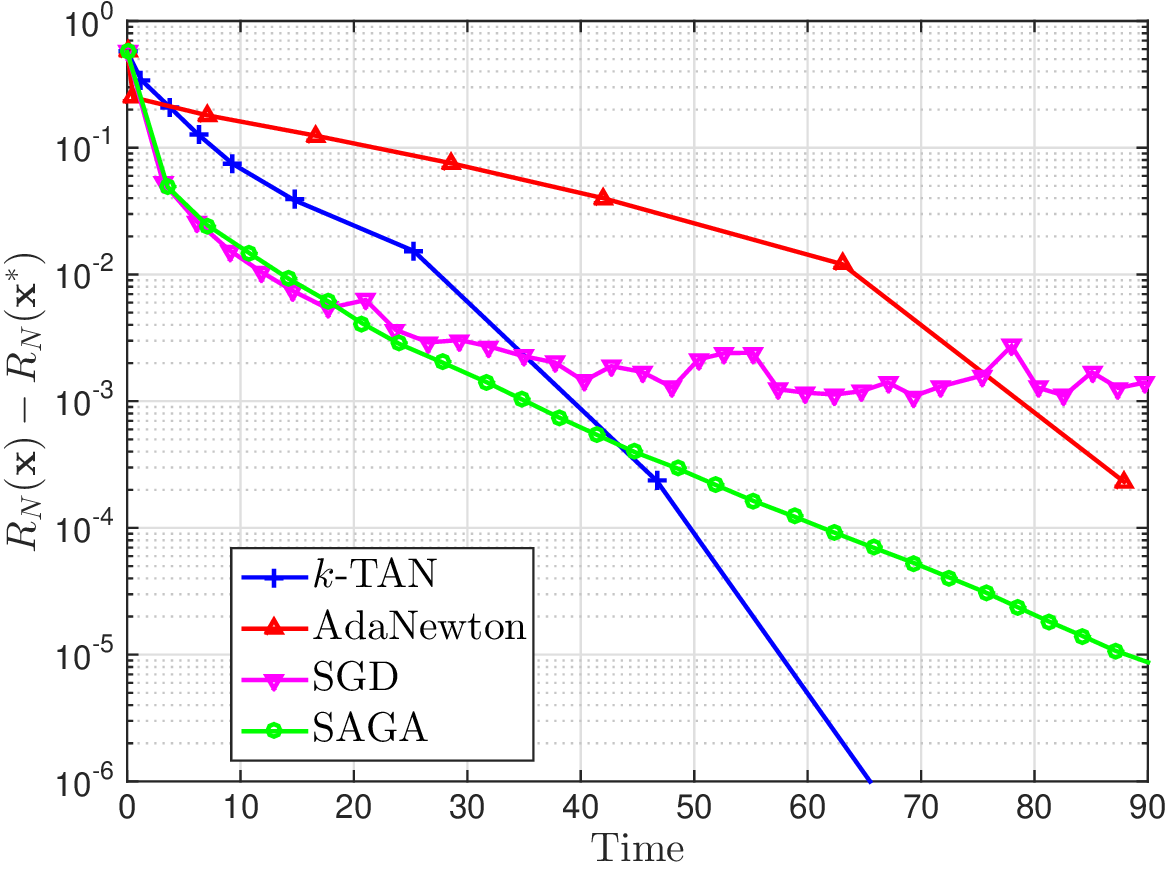}  
          \caption{Convergence of $k$-TAN, AdaNewton, SGD, and SAGA  in terms of number of processed gradients (left) and runtime (right) for the GISETTE handwritten digit classification problem.}
          \label{fig_gisette}
\end{figure}

The GISETTE dataset includes $N=6000$ samples of dimension $p=5000$. We use step sizes of $0.08$ for both SGD and SAGA. In both $k$-TAN and AdaNewton, the sample size is increased by a factor of $\alpha=2$ at each iteration (the condition $\| \nabla R_{n}(\bbw_n)\| < (\sqrt{2 c}) V_n$ is always satisfied) starting with an initial size of $m_0=124$. For both of these methods, we initially run gradient descent on $R_{m_0}$ for $100$ iterations so that we may begin in the statistical accuracy $V_{m_0}$. For $k$-TAN, the truncation $k$ is observed to be able to afford a cutoff of around $0.01p$ in all of our simulations.

In Figure \ref{fig_gisette}, the convergence results of the four methods for GISETTE data is shown. The left plot demonstrates the sub-optimality with respect to the number of gradients, or samples, processed. In particular, $k$-TAN and AdaNewton compute $m$ gradients per iterations, while SGD and SAGA compute 1 gradient per iteration. Observe that the second order methods all converge with a smaller number of total processed gradients than the first order methods, reaching after $2.5\times10^{4}$ samples a sub-optimality of $10^{-7}$. We point out that, while $k$-TAN only approximates the Hessian inverse, its convergence path follows that of AdaNewton exactly. Indeed, both algorithms reach the statistical accuracy of $1/N = 1.6\times10^{-4}$ after 15000 samples, or just over two passes over the dataset. To see the gain in terms of computation time of $k$-TAN over AdaNewton and other methods, we present in the right image of Figure \ref{fig_gisette} the convergence in terms of runtime. In this case, $k$-TAN outperforms all methods, reaching a sub-optimality of $4\times10^{-6}$ after 60 seconds, while AdaNewton reaches only a sub-optimality of $10^{-3}$ after 80 seconds. Note that first order methods have lower cost per iteration than all second order methods. Thus, SAGA is able to converge to $2\times10^{-5}$ after 80 seconds.

 \begin{figure}
	\centering
          \includegraphics[width=0.45\textwidth]{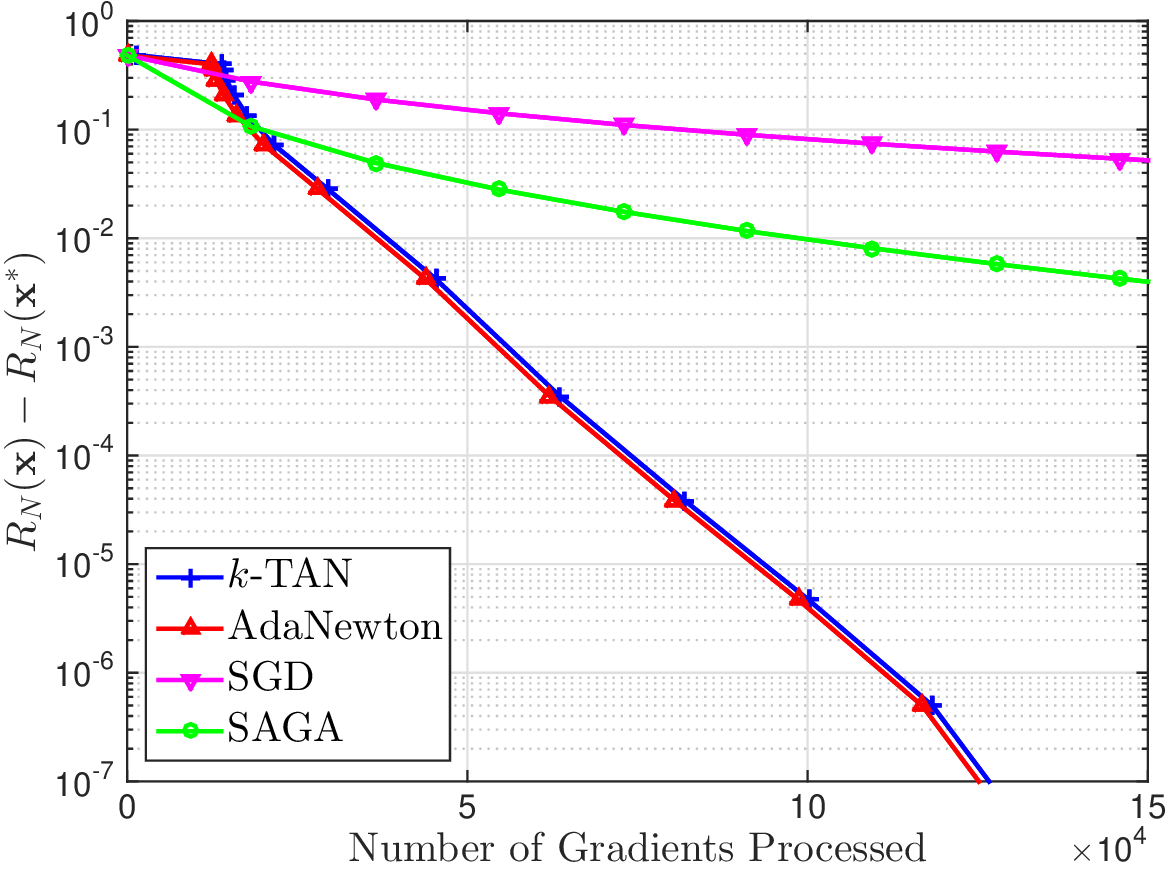}\qquad 
            \includegraphics[width=0.45\textwidth]{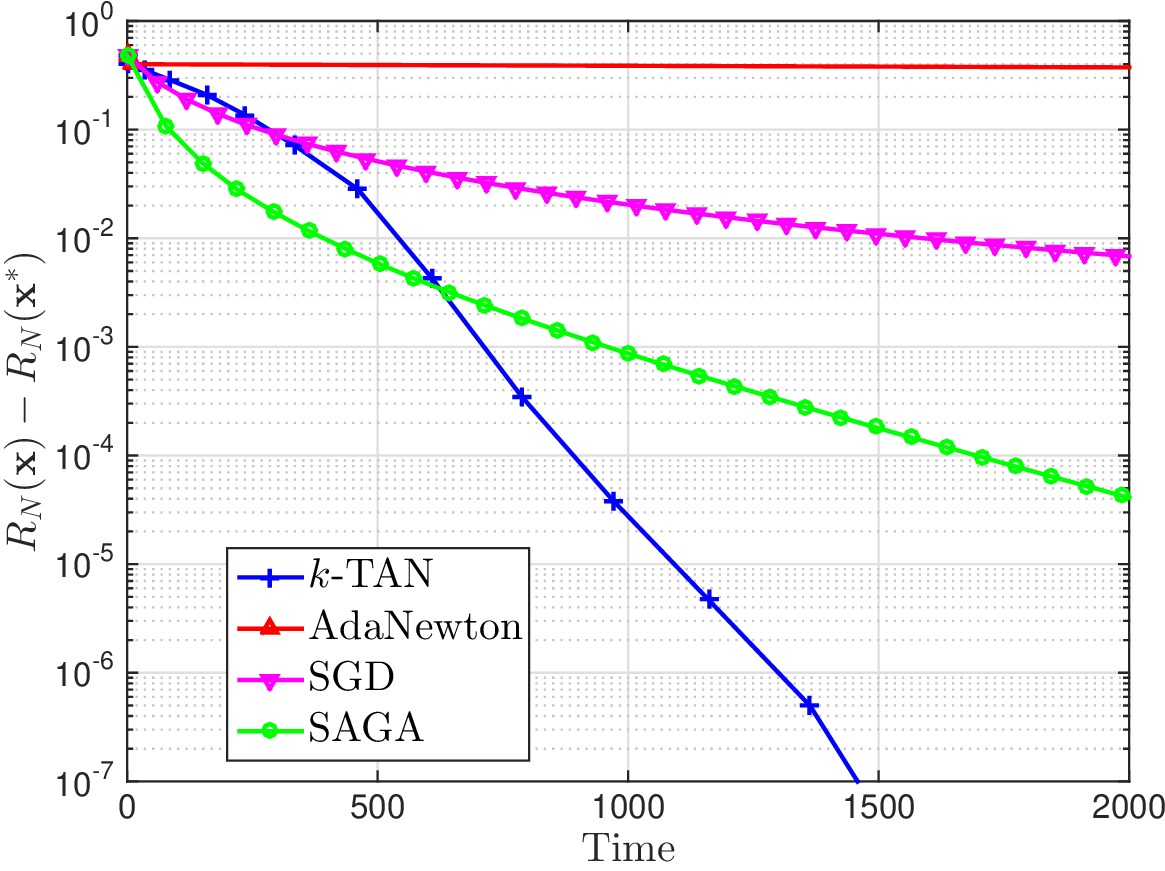}  
          \caption{Convergence of $k$-TAN, AdaNewton, SGD, and SAGA  in terms of number of processed gradients (left) and runtime (right) for the RCV1 text classification problem.}
          \label{fig_rcv1}
\end{figure}

For a very high dimensional problem, we consider the RCV1 dataset with $N=18242$ and $p=47236$. We use step sizes of $0.1$ for both SGD and SAGA and truncate sizes of around $0.001p$ for $k$-TAN, while keeping the parameters for the other methods the same. The results of these simulations are shown in Figure \ref{fig_rcv1}. In the left image, observe that, in terms of processed gradients, the second order methods again outperform the SGD and SAGA, as expected, with $k$-TAN again following the path of AdaNewton. Given the high dimension $p$, the cost of computing the inverse in AdaNewton provides a large bottleneck. The gain in terms of computation time can then be best seen in the right image of Figure \ref{fig_rcv1}. Observe that AdaNewton becomes entirely ineffective in this high dimension. The $k$-TAN method, alternatively, continues to descend at a fast rate because of the inverse truncation step. For this set $k$-TAN outperforms all the other methods, reaching an error of $10^{-7}$ after 1500 seconds. Since both $n$ and $p$ are large, SAGA performs well on this dataset due to small cost per iteration. 

 \begin{figure}
	\centering
          \includegraphics[width=0.45\textwidth]{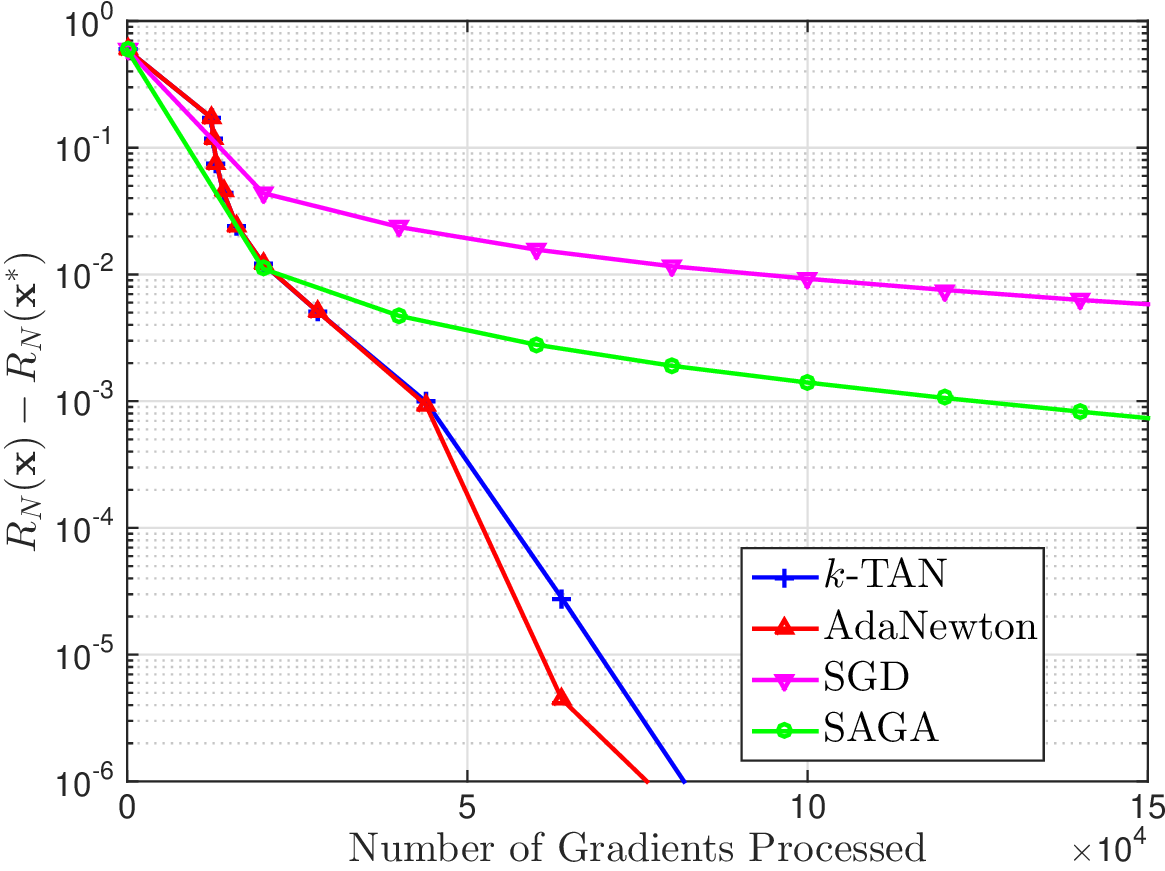}\qquad 
            \includegraphics[width=0.45\textwidth]{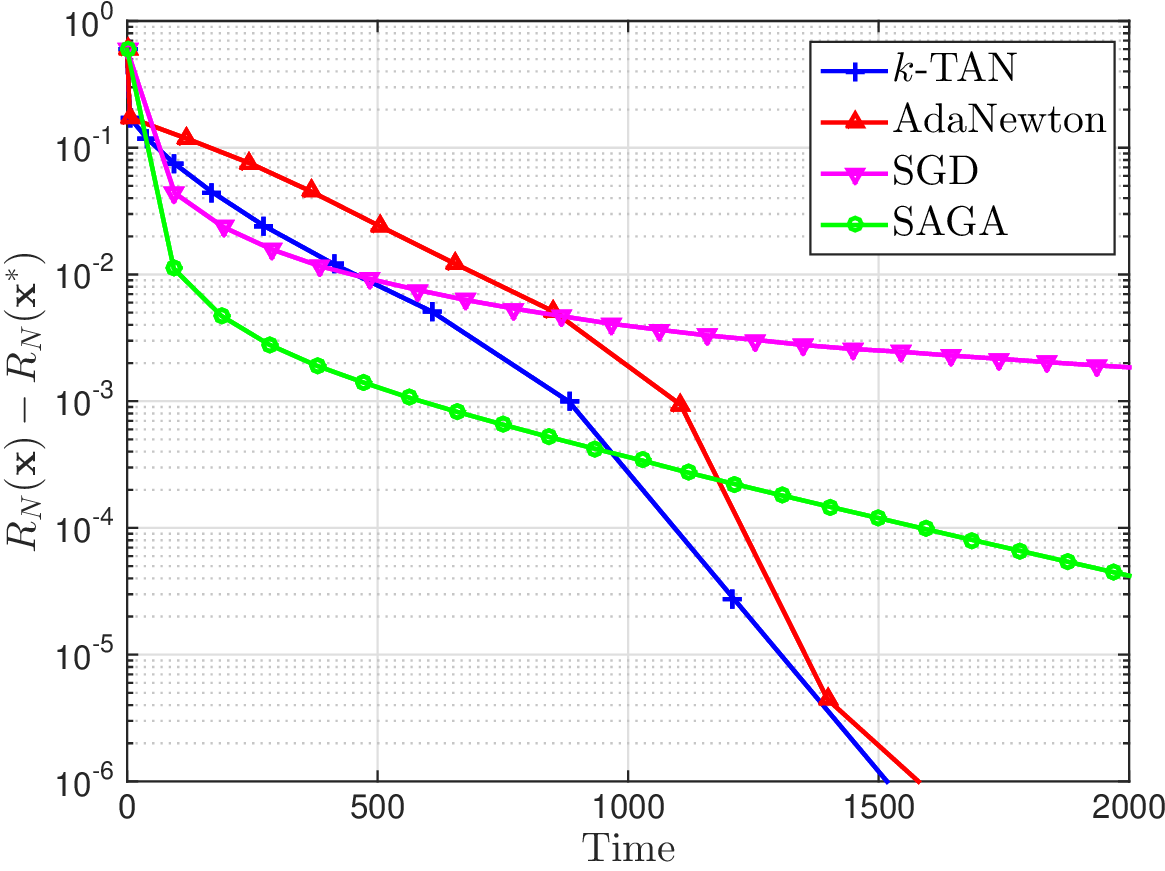}  
          \caption{Convergence of $k$-TAN, AdaNewton, SGD, and SAGA  in terms of number of processed gradients (left) and runtime (right) for the ORANGE text classification problem.}
          \label{fig_orange}
\end{figure}

We perform additional numerical experiments on the ORANGE dataset used for customer relationship prediction in KDD Cup 2009. We use $N=20000$ samples with dimension $p=14472$. The convergence results are shown in Figure \ref{fig_orange}. Observe in the right hand plot that all second order methods, perform similarly well on this dataset. The first order methods, including SAGA, do not converge after 2000 seconds. We also note that, in this experiment, we were able to reduce the truncation size $k$ to around $0.1\%$ of $p$.

 \begin{figure}
	\centering
          \includegraphics[width=0.45\textwidth]{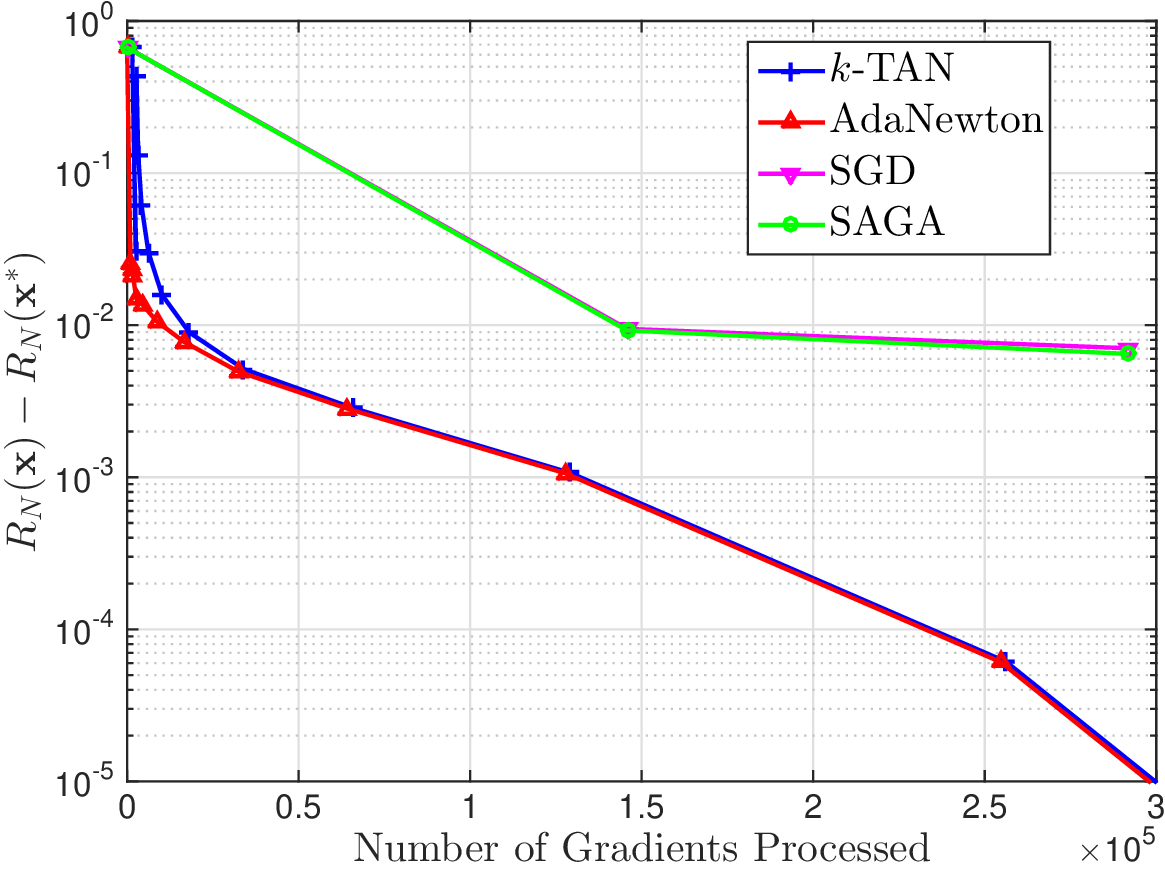}\qquad 
            \includegraphics[width=0.45\textwidth]{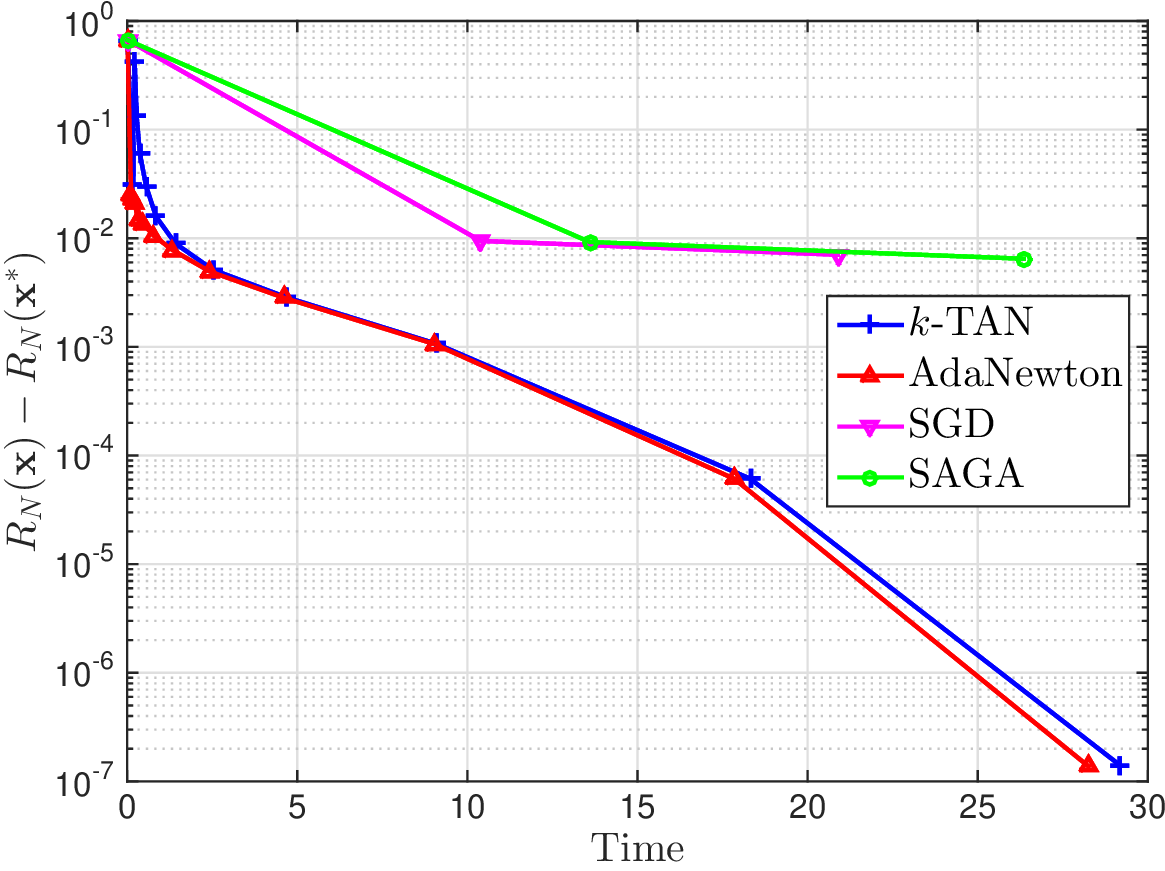}  
          \caption{Convergence of $k$-TAN, AdaNewton, SGD, and SAGA  in terms of number of processed gradients (left) and runtime (right) for the BIO protein homology classification problem.}
          \label{fig_bio}
\end{figure}
 In Figure \ref{fig_bio}, we show results on the BIO dataset used for protein homology classification in KDD Cup 2004. The dimensions are $N=145751$ and $p=74$. In this setting, the number of samples is very large put the problem dimension is very small. Observe in Figure \ref{fig_bio} that both $k$-TAN and AdaNewton greatly outperform the first order methods, due to the reduced cost in Hessian computation that comes from adaptive sample size. However, because $p$ is small, the additional gain from the truncating in the inverse in $k$-TAN does not provide significant benefit relative to AdaNewton.

\section{Discussion}

We demonstrated in this paper the success of the proposed $k$-TAN method on solving large scale empirical risk minimization problems both theoretically and empirically. The $k$-TAN method reduces the total cost in solving \eqref{eq_orig_problem} to its statistical accuracy in two ways: (i) progressively increasing the sample size to reduce the costs of computing gradients and Hessians, and (ii) using a low rank approximation of the Hessian to reduce the cost of inversion. The gain provided by $k$-TAN relative to existing methods is therefore most significant in large scale ERM problems with both large sample size $N$ and dimension $p$. To see this, consider the alternatives previously considered
\begin{itemize}
\item Stochastic first order methods, such as SAGA \citep{defazio2014saga} and SVRG \citep{johnson2013accelerating}, compute a single gradient per iteration, and they have the overall complexity of $\ccalO(N\log(N)p)$ to achieve statistical accuracy of the full training set if $V_N=\mathcal{O}(1/N)$.
\item Newton's method computes gradients and Hessians over the entire dataset and inverts a matrix of size $p$ at each iteration, requiring a total cost of $\ccalO(M (Np^2 + p^3 ))$, where $M$ is number of iterations required to converge. Because Newton's method converges quadratically, $M$ may be small, but the total cost is made large by $Np^2$ and $p^3$.
\item The AdaNewton method \citep{mokhtari2016adaptive} computes gradients and Hessians for a subset of the dataset and inverts a matrix of size $p$ at each iteration, with the size of the subset increasing geometrically. By doubling the sample size every iteration, the statistical accuracy can be reached in a total of $\log_2 N$ iterations after a total of $2N$ passes over the dataset, for a total cost of $\ccalO(2Np^2 + \log_2(N) p^3)$. While Hessian computation cost is reduced, for high dimensional problems the inversion cost of $p^3$ dominates and the algorithm remains costly. 
\end{itemize}

The $k$-TAN method computes gradients and Hessians on a increasing subset of data in the same manner as AdaNewton, put reduces the inversion cost at each iteration to $\ccalO(p^2 \log k )$, resulting in a total cost of $\ccalO(2Np^2 + p^2 \log_2 N \log k )$, or an effective cost of $\ccalO(N p^2)$, if the size of the initial training set is large enough. For ill-conditioned problems, this method is a more feasible option as a second-order method than Newton's method or AdaNewton. This theoretical intuition is indeed supported in the empirical simulations performed on large, high dimensional datasets in this paper.

\acks{We acknowledge the support of the National Science Foundation (NSF CAREER
CCF-0952867) and the Office of Naval Research (ONR N00014-12-1-0997).}

\newpage

\section{Appendix}

\subsection{Proof of Lemma \ref{lemma_sub_bound}}\label{sec_lemma_sub_bound}
We factorize and bound $\| \hbH_{n,k}^{-1}\nabla R_n(\bbx_m) - \bbH_n^{-1}\nabla R_n(\bbx_m) \|$ as
\begin{align} 
\| \hbH_{n,k}^{-1}\nabla R_n(\bbx_m) - \bbH_n^{-1}\nabla R_n(\bbx_m) \| &\leq \| \bbI - \hbH_{n,k}^{-1}\bbH_n \| \|\bbH_n^{-1}\nabla R_n(\bbx_m) \|.
\end{align}
Thus, it remains to bound $\| \bbI - \hbH_{n,k}^{-1}\bbH_n \|$ by some $\epsilon_n$. To do so, consider that we can factorize $\bbH_n = \bbU( \bbSigma + cV_n \bbI)\bbU^T$ and $\hbH^{-1}_n$ as in \eqref{eq_a_hessian_inverse}. We can then expand $\| \bbI - \hbH_{n,k}^{-1}\bbH_n \|$ as
\begin{align}
\| \bbI - \hbH_{n,k}^{-1}\bbH_n \| &= \| \bbI - \bbU [ (\hbSigma_k+cV_n \bbI)^{-1}\times(\bbSigma + cV_n\bbI)]\bbU^T \|,
\end{align}
where $\hbSigma_k \in \reals^{p \times p}$ is the truncated eigenvalue matrix $\bbSigma_k$ with zeros padded for the last $p-k$ diagonal entries. Observe that the first $k$ entries of the product $ (\hbSigma_k+cV_n \bbI)^{-1}\times(\bbSigma + cV_n\bbI)$ are equal to $1$, while the last $p-k$ entries are equal to $(\mu_j + cV_n)/cV_n$. Thus, we have that 
\begin{align}
\| \bbI - \hbH_{n,k}^{-1}\bbH_n \| &= \left| \frac{\mu_{k+1}}{cV_n} \right|.
\end{align}

\subsection{Proof of Lemma \ref{lemma_sub_n}}

To begin, recall the result from Lemma \ref{lemma_sub_bound} in \eqref{eq_lemma_sub_bound}. From this, we use the following result from \citep[Lemma 6]{pilanci2015newton}, which present here as a lemma.

\begin{lemma}\label{lemma_decr_sub}
Consider the $k$-TAN step where $\| \hbH_{n,k}^{-1}\nabla R_n(\bbx_m) - \bbH_n^{-1} \nabla R_n(\bbx_m)\| \leq \epsilon_n \| \bbH_n^{-1} \nabla R_n(\bbx_m)\| $. The Newton decrement of the $k$-TAN iterate $\lambda_n(\bbx_n)$ is bounded by
\begin{align}\label{eq_lemma_decr_sub}
\lambda_n(\bbx_n) \leq \frac{1}{(1 - (1+\epsilon_n) \lambda_n(\bbx_m))^2}\left[ (1+\epsilon_n) \lambda_n(\bbx_m)^2 + \epsilon_n \lambda_n(\bbx_m) \right] \quad w.h.p
\end{align}
\end{lemma}

Lemma \ref{lemma_decr_sub} provides a bound on the Newton decrement of the iterate $\bbx_n$ computed from the $k$-TAN update in \eqref{eq_update} in terms of Newton decrement of the previous iterate $\bbx_m$ and the error $\epsilon_n$ incurred from the truncation of the Hessian. We proceed in a manner similar to \citep[Proposition 4]{ mokhtari2016adaptive} by finding upper and lower bounds for the sub-optimality $S_n(\bbx) = R_n(\bbx)-R_n(\bbx_n^*)$ in terms of the Newton decrement parameter $\lambda_n (\bbx)$. Consider the result from \citep[Theorem 4.1.11]{nesterov1998introductory},
 \begin{equation}\label{proof_of_prop_10}
		  \lambda_n (\bbx) - \ln\left(1+\lambda_n (\bbx)\right) \leq  
	 R_n(\bbx)-R_n(\bbx_n^*)  \leq
	  - \lambda_n (\bbx) - \ln\left(1-\lambda_n (\bbx)\right).
\end{equation}

Consider the Taylor's expansion of $\ln(1+a)$ for $a=\lambda_n (\bbx)$ to obtain the lower bound on $\lambda_n (\bbx)$, 
\begin{equation}\label{eq_taylor}
\lambda_n (\bbx) \geq \ln\left(1+\lambda_n (\bbx)\right) + \frac{1}{2}\lambda_n (\bbx)^2-\frac{1}{3}\lambda_n (\bbx)^3.
\end{equation}
Assume that $\bbx$ is such that $0<\lambda_n (\bbx)<1/4$. Then the expression in \eqref{eq_taylor} can be rearranged and bounded as
\begin{equation}\label{eq_taylor2}
\frac{1}{6} \lambda_n (\bbx)^2 \leq \frac{1}{2}\lambda_n (\bbx)^2- \frac{1}{3}\lambda_n (\bbx)^3
\end{equation}
Now, consider the Taylor's expansion of $\ln(1-a)$ for  $a=\lambda_n (\bbx)$ in a similar manner to obtain for $\lambda_n (\bbx)<1/4$, from \citep[Chapter 9.6.3]{boyd04}.
\begin{equation}\label{eq_taylor3}
- \lambda_n (\bbx) - \ln\left(1-\lambda_n (\bbx)\right) \leq \lambda_n (\bbx)^2
\end{equation}
Using these bounds with the inequalities in \eqref{proof_of_prop_10} we obtain the upper and lower bounds on $S_n(\bbx)$ as  
\begin{align}\label{eq_subopt_decr}
\frac{1}{6}\lambda_n(\bbx)^2 \leq S_n(\bbx) \leq \lambda_n(\bbx)^2.
\end{align}
Now, consider the bound for Newton decrement of the $k$-TAN iterate $\lambda_n(\bbx_n)$ from \eqref{eq_lemma_decr_sub}. As we assume that $\lambda_n (\bbx_m) < 1/4$, we have
\begin{align}
\lambda_n(\bbx_n) \leq \frac{4}{(3 - \epsilon_n)^2}\left[ (1+\epsilon_n) \lambda_n(\bbx_m)^2 + \lambda_n(\bbx_m)\epsilon_n\right].
\end{align}
We substitute this back into the upper bound in \eqref{eq_subopt_decr} for $\bbx = \bbx_n$ to obtain
\begin{align}\label{eq_subopt_decr1}
S_n(\bbx_n) &\leq \lambda_n(\bbx_n)^2 \leq \frac{16}{(3 - \epsilon_n)^4}\left[ (1+\epsilon_n) \lambda_n(\bbx_m)^2 + \lambda_n(\bbx_m)\epsilon_n \right]^2 \nonumber \\
& = \frac{16}{(3 - \epsilon_n)^4}\left[ (1+\epsilon_n)^2 \lambda_n(\bbx_m)^4 +2 \epsilon_n(1+\epsilon_n) \lambda_n(\bbx_m)^3 +  \epsilon_n^2 \lambda_n(\bbx_m)^2 \right]. 
\end{align}
Consider also from \eqref{eq_subopt_decr} that we can upper bound the Newton decrement as $\lambda(\bbx_m)^2 \leq 6 S_n(\bbx_m)$. We plug this back into \eqref{eq_subopt_decr1} to obtain a final bound for sub-optimality as
\begin{align}\label{eq_subopt_decr2}
S_n(\bbx_n) \leq \frac{16}{(3 - \epsilon_n)^4}\left[ 36(1+\epsilon_n)^2 S_n(\bbx_m)^2 + 30 \epsilon_n(1+\epsilon_n)S_n(\bbx_m)^{3/2} + 6\epsilon_n^2 S_n(\bbx_m) \right]. 
\end{align}

\subsection{Proof of Theorem \ref{theorem_main_result}}\label{sec_theorem_main_result}

The proof of this theorem follows from the previous results. Observe that, from Proposition \ref{main_prop} the condition in \eqref{cond_1} ensures that $\bbx_m$ will be in the quadratic region of $R_n$, i.e. $\lambda_n(\bbx_m) < 1/4$. This condition validates the result in \eqref{eq_lemma_sub_n}, restated as
\begin{align}\label{eq_subopt_decr2_a}
S_n(\bbx_n) \leq \frac{16}{(3 - \epsilon_n)^4}\left[ 36(1+\epsilon_n)^2 S_n(\bbx_m)^2 + 30 \epsilon_n(1+\epsilon_n)S_n(\bbx_m)^{3/2} + 6\epsilon_n^2 S_n(\bbx_m) \right]
\end{align}

From Proposition \ref{prop_sub_m} we can bound the $n$-suboptimality of the previous iterate $S_n(\bbx_m)$. For notational simplicity, we focus on the case in which the statistical accuracy is  $V_m = \ccalO(1/m)$, as given in \eqref{eq_prop_sub_m_2}. Furthermore, we can take the expression given for the truncation error $\epsilon_n$ from \eqref{eq_lemma_sub_bound}. Consider for some $\rho$, the $k+1$-th eigenvalue of the Hessian satisfies $\mu_{k+1} \leq \rho cV_n$. Substituting this and \eqref{eq_prop_sub_m_2} into \eqref{eq_subopt_decr2_a}, we obtain
\begin{align}\label{eq_subopt_new}
S_n(\bbx_n) \leq \frac{16}{(3 - \rho)^4}\left[ 36 K^2(1+\rho)^2 V_m^2 + 30 K^{3/2} \rho(1+\rho)V_m^{3/2} + 6K\rho^2 V_m\right].
\end{align}
The bound in \eqref{eq_subopt_new} provides us then the condition in\eqref{cond_2} for $S_n(\bbx_n) \leq V_n$.

 \newpage
 
 \bibliography{bibliography}
\bibliographystyle{icml2014}
\end{document}